\documentclass[a4paper,12pt]{article}

\usepackage{authblk}
\usepackage{graphicx}
\usepackage[truedimen,margin=2.5truecm]{geometry}
\usepackage{amsthm,amssymb,amsmath,mathrsfs,xcolor}

\numberwithin{equation}{section} 
\newtheorem{theorem}{Theorem}[section]
\newtheorem{corollary}[theorem]{Corollary}

\newtheorem{proposition}[theorem]{Proposition}

\theoremstyle{definition}
\newtheorem{remark}[theorem]{Remark}

\newcommand{\C}{\mathbb{C}} 
\newcommand{\R}{\mathbb{R}} 
\newcommand{\N}{\mathbb{N}} 

\title{Topology and convergence on the space of measure-valued functions}
\author{Takahiro Hasebe}
\affil{\small Department of Mathematics, Hokkaido University, North 10 West 8, Kita-ku, Sapporo 060-0810, Japan.
\texttt{thasebe@math.sci.hokudai.ac.jp}}
\author{Ikkei Hotta}
\affil{\small Department of Applied Science, Yamaguchi University 2-16-1 Tokiwadai, Ube 755-8611, Japan.
\texttt{ihotta@yamaguchi-u.ac.jp}}
\author{Takuya Murayama}
\affil{\small Faculty of Mathematics, Kyushu University, 744 Motooka, Nishi-ku, Fukuoka 819-0395, Japan ({\it After completing this manuscript, Murayama moved to} Department of Mathematics, Graduate School of Science, Kobe University, 1-1 Rokkodai, Nada-ku, Kobe 657-8501, Japan.)
\texttt{murayama@math.kobe-u.ac.jp}}
\date{}

\begin{document}

\maketitle

\begin{abstract}
In these notes, uniform convergence on compacta is studied on the space of functions taking values in the set of finite Borel measures.
Related limit theorems, including L\'evy's continuity theorem and functional limit theorems for (classical and non-commutative) additive processes, are also described.
N.B.: the contents of this manuscript have been incorporated
into another manuscript (arXiv:2412.18742). 


\end{abstract}

\section{Introduction}
\label{sec:introduction}

In the theory of stochastic processes, a family of measures $\mu_t$ indexed by a parameter $t$, such as time or initial point, naturally appears. Such a family can be regarded as a measure-valued function $t \mapsto \mu_t$ on the parameter set.  This paper offers some simple and fundamental observations on the topology and convergence on the space of measure-valued continuous functions. 

Let $S,T$ be topological spaces. Let $\mathbf{M}(S)$ denote the set of nonnegative finite Borel measures   on $S$ equipped with the weak topology. We designate the set of continuous mappings from $T$ to the space $\mathbf{M}(S)$ as $C(T; \mathbf{M}(S))$.
Our aim is to find a topology on $C(T; \mathbf{M}(S))$ such that the sequential convergence $(\mu^{(n)}_t)_{t \in T} \to (\mu_t)_{t \in T}$ ($n \to \infty$) in that topology occurs if and only if for every compact set $K \subset T$ and for every bounded continuous function $f \colon S \to \R$ we have
\begin{equation} \label{eq:seq-luwc}
\lim_{n \to \infty}\sup_{t \in K} \left\lvert \int_S f(x) \, \mu^{(n)}_t(dx) - \int_S f(x) \, \mu_t(dx) \right\rvert =0.
\end{equation}
The expression \eqref{eq:seq-luwc} can be regarded as the ``uniform convergence of $((\mu^{(n)}_t)_{t \in T})_{n \in \N}$ on the compact set $K$.''
However, since a topology is in general not uniquely determined by sequential convergence in that topology, there might be more than one topologies satisfying our requirement.
To be precise, therefore, the sequence $((\mu^{(n)}_t)_{t \in T})_{n \in \N}$ above should be replaced by an arbitrary \emph{net} $((\mu^{(\alpha)}_t)_{t \in T})_{\alpha \in A}$.  After introducing a topology on $C(T; \mathbf{M}(S))$ in Section \ref{sec:cpt-open_top},  we consider relationships between the topology and convergence of characteristic functions in Section \ref{sec:cf}, and then vague convergence of probability measures and convergence of moments in Section \ref{sec:vague_moment}. 
Sections~\ref{sec:FLT} and \ref{sec:BP_bij} are devoted to applications, which originally motivated the present work.  The convergence \eqref{eq:seq-luwc} naturally appears in a functional limit theorem for additive processes, see Section~\ref{sec:FLT}, and also in the convergence  of the marginal distributions of noncommutative stochastic processes, see Section~\ref{sec:BP_bij}.

\section{Topology on $C(T; \mathbf{M}(S))$}
\label{sec:cpt-open_top}

In what follows, we regard $\mathbf{M}(S)$ and $C(T; \mathbf{M}(S))$ as \emph{uniform spaces} to obtain the topology $\mathscr{T}$ prescribed above instead of defining it directly via nets.
The reader can consult Kelley~\cite[Chapters 6 and 7]{Ke75} for uniform spaces and uniform convergence of functions.
We also see that $\mathscr{T}$ is metrizable and hence ``determined by sequences'' if $S$ is a separable metric space (Corollary~\ref{cor:metrizability}).

Let $C_{\rm b}(S)$ denote the set of bounded continuous  real-valued functions on $S$.
Define $\mathbf{U}$ as the uniformity (or uniform structure) on $\mathbf{M}(S)$ generated by the collection of sets
\begin{equation} \label{eq:equi-distant_set}
U(f, \varepsilon)
:=\left\{\, (\mu, \nu)\in \mathbf{M}(S)^2 : \left\lvert \int f \,d\mu - \int f \,d\nu \right\rvert<\varepsilon \,\right\}
\end{equation}
indexed by $f \in C_{\rm b}(S)$ and $\varepsilon>0$.
The topology of the uniform space $(\mathbf{M}(S), \mathbf{U})$ is exactly the weak topology.
We further define $\mathscr{U}$ as the uniformity on $C(T; \mathbf{M}(S))$ generated by the collection of sets
\begin{align}
V(K; f, \varepsilon)
&:=\left\{\, (\mu_t, \nu_t)_{t\in T}\in C(T; \mathbf{M}(S))^2
: (\mu_t, \nu_t)\in U(f, \varepsilon)\ \text{for all}\ t\in K\,\right\}
\notag \\
&=\left\{\, (\mu_t, \nu_t)_{t\in T}\in C(T; \mathbf{M}(S))^2
: \sup_{t\in K}\left\lvert \int f \,d\mu_t - \int f \,d\nu_t \right\rvert<\varepsilon \,\right\}
\label{eq:lu_equi-distant_set}
\end{align}
indexed by compact $K \subset T$, $f \in C_{\rm b}(S)$, and $\varepsilon>0$.
Then the topology $\mathscr{T}$ on $C(T; \mathbf{M}(S))$ induced by $\mathscr{U}$ is, by definition, the topology of ($\mathbf{U}$-)uniform convergence on compacta%
\footnote{``Uniform convergence on compacta'' means the same as ``locally uniform convergence'' if $T$ is locally compact.}.
Obviously, $\lim_\alpha (\mu^{(\alpha)}_t)_{t \in T}=(\mu_t)_{t \in T}$%
\footnote{In the general setting here, the limit $(\mu_t)_{t \in T}$ may not be unique because no separation axiom is assumed.
However, this is not our focus, and actually, $S$ will be assumed to be metrizable later.}
in $\mathscr{T}$ if and only if
\begin{equation} \label{eq:net-luwc}
\lim_\alpha\sup_{t \in K} \left\lvert \int_S f(x) \, \mu^{(\alpha)}_t(dx) - \int_S f(x) \, \mu_t(dx) \right\rvert=0
\end{equation}
for every compact set $K \subset T$ and for every $f \in C_{\rm b}(S)$.

By the following fact, adopted from Kelley~\cite[Theorem~11, Chapter~7]{Ke75}, $\mathscr{T}$ is indeed the \emph{compact-open topology} \cite[p.221]{Ke75} on $C(T; \mathbf{M}(S))$:

\begin{proposition} \label{prop:cpt-open}
Let $X$ be a topological space and $(Y, \mathcal{U})$ be a uniform space.
On the set $C(X,Y)$ of continuous mappings from $X$ to $Y$, the topology of $\mathcal{U}$-uniform convergence on compacta is identical with the compact-open topology.
In particular, the former topology depends only on the topology that $\mathcal{U}$ induces on $Y$.
\end{proposition}

Proposition~\ref{prop:cpt-open} enables us to obtain the topology $\mathscr{T}$ in a way different from above.
For example, let us suppose that $S$ is a separable metric space.
Then there exists a distance $\rho$ on $\mathbf{M}(S)$ which induces the weak topology (e.g., the L\'evy--Prokhorov distance).
We replace $U(f, \varepsilon)$ in \eqref{eq:equi-distant_set} with
\[
U_\rho(\varepsilon):=\{\, (\mu, \nu) \in \mathbf{M}(S)^2 : \rho(\mu, \nu)<\varepsilon \,\}
\]
and define $V_\rho(K; \varepsilon)$ in a way similar to \eqref{eq:lu_equi-distant_set}, which leads to another uniformity $\mathscr{U}_\rho$ on $C(T; \mathbf{M}(S))$.
By Proposition~\ref{prop:cpt-open}, $\mathscr{U}_\rho$ induces the topology $\mathscr{T}$.
Hence we have the following:

\begin{corollary} \label{cor:metrizability}
Let $S$ be a separable metric space and $T$ be a topological space.
\begin{enumerate}
\item \label{i:lu_LP}
Fix a distance $\rho$ which induces the weak topology on $\mathbf{M}(S)$.
For a net $((\mu^{(\alpha)}_t)_{t \in T})_{\alpha \in A}$ in $C(T; \mathbf{M}(S))$ and $(\mu_t)_{t \in T}\in C(T; \mathbf{M}(S))$, the following are equivalent:
\begin{itemize}
\item
\eqref{eq:net-luwc} holds for every compact $K \subset T$ and for every $f \in C_{\rm b}(S)$;
\item
$\lim_\alpha\sup_{t \in K}\rho(\mu^{(\alpha)}_t, \mu_t)=0$ holds for every compact $K$.
\end{itemize}

\item \label{i:luwc_metrizability}
Suppose that $T$ is locally compact and $\sigma$-compact.
Then the topology $\mathscr{T}$ on $C(T; \mathbf{M}(S))$ is metrizable.
\end{enumerate}
\end{corollary}

\begin{proof}
\eqref{i:lu_LP} merely states that $\mathscr{U}$ and $\mathscr{U}_\rho$ induce the same topology.
In \eqref{i:luwc_metrizability} a desired distance is given by
\[
\mathscr{D}_\rho\left((\mu_t)_t, (\nu_t)_t \right):=\sum_{j\in \mathbb{N}}2^{-j}\left( \sup_{t\in K_j}\rho(\mu_t, \nu_t)\wedge 1\right),
\qquad (\mu_t)_{t \in T}, (\nu_t)_{t \in T} \in C(T;\mathbf{M}(S)),
\]
as usual.
Here, $(K_j)_{j \in \N}$ is an exhaustion sequence of $T$.
\end{proof}

\section{Characteristic function and L\'evy's continuity theorem}
\label{sec:cf}

If $S$ is the $d$-dimensional Euclidean space $\R^d$, we can obtain the topology $\mathscr{T}$ in Section~\ref{sec:cpt-open_top} by means of characteristic functions.
Writing the characteristic function of $\mu \in \mathbf{M}(\R^d)$ as $\hat{\mu}(\xi):=\int e^{i\xi \cdot x} \,\mu(dx)$, we identify the space $\mathbf{M}(\R^d)$ with the space $\widehat{\bf M}(\R^d):=\{\, \hat{\mu} : \mu \in \mathbf{M}(\R^d) \,\}$ of characteristic functions.
The latter space is a subset of $C(\R^d; \C)$ and hence has a distance $\hat{\rho}$ which induces the topology of locally uniform convergence.
By L\'evy's continuity theorem, the distance $\hat{\rho}$, transferred to $\mathbf{M}(\R^d)$ through our identification, induces the weak topology.
Thus, we can take $\hat{\rho}$ as the distance $\rho$ discussed in Section~\ref{sec:cpt-open_top}.
To put it differently, given a topological space $T$, the homeomorphism $\mathbf{M}(\R^d) \ni \mu \mapsto \hat{\mu} \in \widehat{\bf M}(\R^d)$ is lifted to the homeomorphism $C(T; \mathbf{M}(\R^d)) \ni (\mu_t)_t \mapsto (\widehat{\mu_t})_t \in C(T; \widehat{\bf M}(\R^d))$.

We can think that the preceding paragraph generalizes L\'evy's continuity theorem for measure-valued functions.
In this direction, we can prove the following:

\begin{theorem} \label{th:Levy_continuity}
Let $T$ be a first-countable Hausdorff space and $\rho$ be a distance on $\mathbf{M}(\R^d)$ which induces the weak topology.
Suppose that $\mu^{(n)}_t \in \mathbf{M}(\R^d)$ for each $t \in T$ and for each $n \in \N$ and $(\mu_t)_{t \in T} \in C(T; \mathbf{M}(\R^d))$.
Then the following are mutually equivalent:
\begin{enumerate}
\item \label{i:Levy_luwc}
$\lim_{n \to \infty}\sup_{t \in K}\left\lvert \int_{\R^d} f \, d\mu^{(n)}_t - \int_{\R^d} f \, d\mu_t \right\rvert=0$
for every compact set $K \subset T$ and for every $f \in C_{\rm b}(\R^d)$;
\item \label{i:Levy_dist}
$\lim_{n \to \infty}\sup_{t \in K}\rho(\mu^{(n)}_t, \mu_t)=0$
for every compact set $K \subset T$;
\item \label{i:Levy_cf-lu}
$(\mu^{(n)}_t)^\land(\xi) \to \widehat{\mu_t}(\xi)$ ($n \to \infty$) uniformly in $(\xi, t)$ on each compact subset of $\mathbb{R}^d \times T$;
\item \label{i:Levy_cf-pw}
for each fixed $\xi \in \mathbb{R}^d$, $(\mu^{(n)}_t)^\land(\xi) \to \widehat{\mu_t}(\xi)$ ($n \to \infty$) uniformly in $t$ on each compact subset of $T$.
\end{enumerate}
\end{theorem}

\begin{proof}
We begin with a simple observation.
Let $K \subset T$ be compact and $(t(n))_{n \in \N}$ be a sequence in $K$ which converges to some $u$.
By assumption, $\mu_{t(n)} \xrightarrow{\rm w} \mu_u$ as $n \to \infty$.
Then each one of \eqref{i:Levy_luwc}--\eqref{i:Levy_cf-pw} yields $\mu^{(n)}_{t(n)} \xrightarrow{\rm w} \mu_u$.
For example, if \eqref{i:Levy_cf-pw} holds, then for each $\xi \in \R^d$,
\begin{align*}
\lvert (\mu^{(n)}_{t(n)})^\land(\xi) - \widehat{\mu_u}(\xi) \vert
&\le \lvert (\mu^{(n)}_{t(n)})^\land(\xi) - \widehat{\mu_{t(n)}}(\xi) \vert + \lvert \widehat{\mu_{t(n)}}(\xi) - \widehat{\mu_u}(\xi) \vert \\
&\le \sup_{t \in K}\lvert (\mu^{(n)}_t)^\land(\xi) - \widehat{\mu_t}(\xi) \rvert + \lvert \widehat{\mu_{t(n)}}(\xi) - \widehat{\mu_u}(\xi) \vert
\to 0,
\end{align*}
which implies $\mu^{(n)}_{t(n)} \xrightarrow{\rm w} \mu_u$ by L\'evy's continuity theorem.

Now assume that one of \eqref{i:Levy_luwc}--\eqref{i:Levy_cf-pw} fails.
Then there exist a compact set $K \subset T$, a positive number $\varepsilon$, and an increasing sequence $(n(k))_{k \in \N}$ in $\N$ such that $\sup_{t \in K}D(\mu^{(n(k))}_t, \mu_t) \ge 2\varepsilon$.
Here, $D(\mu, \nu)$ is one of
\[
\left\lvert \int f \,d\mu - \int f \,d\nu \right\rvert,\ \rho(\mu, \nu),\ \sup_{\xi \in B}\lvert \hat{\mu}(\xi) - \hat{\nu}(\xi) \rvert,\ \lvert \hat{\mu}(\xi) - \hat{\nu}(\xi) \rvert
\]
for some $f \in C_{\rm b}(S)$ in the first and for some compact $B \subset \R^d$ in the third expression, according to which of \eqref{i:Levy_luwc}--\eqref{i:Levy_cf-pw} is assumed to fail.
Then, taking a further subsequence of $(n(k))_{k \in \N}$ if necessary, we can choose a sequence $(t(k))_{k \in \N}$ in $K$ which converges to some $u$ so that $D(\mu^{(n(k))}_{t(k)}, \mu_{t(k)}) \ge \varepsilon$.
This implies that all of \eqref{i:Levy_luwc}--\eqref{i:Levy_cf-pw} fail by the previous observation.
Indeed, if $\mu_{t(k)} \xrightarrow{\rm w} \mu_u$ and $\mu^{(n(k))}_{t(k)} \xrightarrow{\rm w} \mu_u$ occur at the same time, then $D(\mu^{(n(k))}_{t(k)}, \mu_{t(k)})$ must go to zero as $k \to \infty$.
Thus the theorem is proved.
\end{proof}

\begin{remark} \label{rem:top_seq}
By the first paragraph of this section, Conditions \eqref{i:Levy_luwc}--\eqref{i:Levy_cf-lu} are mutually equivalent if $((\mu^{(n)}_t)_{t \in T})_{n \in \N}$ is replaced by a net $((\mu^{(\alpha)}_t)_{t \in T})_{\alpha \in A}$ in $C(T; \mathbf{M}(S))$ as in Corollary~\ref{cor:metrizability}.
In Theorem~\ref{th:Levy_continuity}, we have been concerned with sequences only so as to remove the $t$-continuity assumption on $((\mu^{(n)}_t)_{t \in T})_{n \in \N}$ and, furthermore, to include the $\xi$-pointwise condition~\eqref{i:Levy_cf-pw}.
Thus, the topological method in Section~\ref{sec:cpt-open_top} and sequential method in the present section do not completely include each other.
\end{remark}

\section{Vague convergence of probability measures and convergence of moments}
\label{sec:vague_moment}

Keeping Remark~\ref{rem:top_seq} in mind, we follow the ideas in Sections~\ref{sec:cpt-open_top} and \ref{sec:cf} to discuss two more situations, in which the values of the functions at issue are probability measures.  The proof of Theorems~\ref{th:weak_vs_vague} and \ref{th:Wasserstein} below is quite similar to that of Theorem~\ref{th:Levy_continuity} and thus omitted.  Let $\mathbf{P}(S)$ be the set of Borel probability measures on a topological space $S$. 

First, we focus on the vague convergence of probability measures.
Suppose that $S$ is a locally compact second countable Hausdorff space.
This is the same as saying that $S$ is a locally compact \emph{Polish} space (i.e., space which is separable and metrizable with a complete metric); e.g., see \cite[Remark~5, {\S}29]{Ba01}.
The weak topology on $\mathbf{M}(S)$ is then metrizable (with the Levy--Prokhorov distance as we have already mentioned).
Now let $C_\infty(S)$ (resp.\ $C_{\rm c}(S)$) denote the set of continuous real functions on $S$ vanishing at infinity (resp.\ with compact support), and define $\mathbf{U}^\prime$ (resp. $\mathbf{U}^{\prime\prime}$) as the uniformity on $\mathbf{M}(S)$ generated by the collection of sets \eqref{eq:equi-distant_set} indexed by $f \in C_\infty(S)$ (resp.\ $C_{\rm c}(S)$) and $\varepsilon>0$.
Then the topology induced by $\mathbf{U}^{\prime\prime}$ is called the vague topology on $\mathbf{M}(S)$.
As is seen in Remark~\ref{rem:vague} below, on $\mathbf{P}(S)$ the vague topology coincides with the weak one; in other words, the uniformities $\mathbf{U}$, $\mathbf{U}^\prime$, and $\mathbf{U}^{\prime\prime}$ restricted to $\mathbf{P}(S)$ induce the same topology.
By Proposition~\ref{prop:cpt-open}, the topology of uniform convergence on compacta on the space $C(T; \mathbf{P}(S))$ is the same if we take any one of these uniformities.

\begin{remark} \label{rem:vague}
A nice exposition on the vague topology can be found in, e.g., Bauer's book~\cite[{\S}30--31]{Ba01}.
As a premise, it is natural to consider the vague topology on the set of \emph{Radon} measures, but in our case, every finite Borel measure on $S$ is indeed Radon because $S$ is assumed to be a locally compact second countable Hausdorff space.
Then it is well known \cite[Corollary~30.9]{Ba01} that a sequence $(\mu_n)_{n \in \N}$ in $\mathbf{P}(S)$ converges vaguely to a measure $\mu \in \mathbf{P}(S)$ if and only if it converges weakly to $\mu$.
This implies that the identity mapping on $\mathbf{P}(S)$, with the domain and codomain equipped with the vague and weak topologies, respectively, is sequentially bi-continuous.
Moreover, the vague topology on $\mathbf{P}(S)$ is metrizable \cite[Theorem~31.5]{Ba01} as well as the weak one.
Thus, the identity mapping above is homeomorphic, which means that the vague topology coincides with the weak one on $\mathbf{P}(S)$.
Note that this is no longer valid on the larger space $\mathbf{M}(S)$; see examples in the same book \cite[Remarks 1 and 2, {\S}30]{Ba01}.
\end{remark}

The sequential counterpart to the preceding topological discussion is the following:

\begin{theorem} \label{th:weak_vs_vague}
Let $S$ be a locally compact second-countable Hausdorff space and $T$ be a first-countable Hausdorff space.
Suppose that $\mu^{(n)}_t \in \mathbf{P}(S)$ for each $t \in T$ and for each $n \in \N$ and $(\mu_t)_{t \in T} \in C(T; \mathbf{P}(S))$.
Designate the following condition as $(C_\star)$ for $\star \in \{{\rm b}, \infty, {\rm c}\}$: \eqref{eq:seq-luwc} holds for every compact $K \subset T$ and for every $f \in C_\star(S)$.
Then the three conditions $(C_\star)$, $\star \in \{{\rm b}, \infty, {\rm c}\}$, are mutually equivalent.
\end{theorem}

Next, we focus on the convergence of moments.
Let $(S, d)$ be a complete separable metric space.
For a fixed $p \ge 1$, we define $\mathbf{P}^p(S, d)$ as the set of Borel probability measures with finite $p$-th moment. More precisely, $\mathbf{P}^p(S, d) := \{\, \mu \in \mathbf{P}(S) : \int_S d(x_0, x)^p \,\mu(dx)<+\infty \,\}$ for some $x_0 \in S$.
Here, the definition does not depend on the choice of $x_0$.
For $\mu \in \mathbf{P}^p(S, d)$, we put $m_p(\mu):=\int_S d(x_0, x)^p \,\mu(dx)$.

On $\mathbf{P}^p(S, d)$ a uniformity $\mathbf{U}^p$ is generated by the collection of sets
\[
U_p(f,\varepsilon)
:=\left\{\, (\mu, \nu)\in \mathbf{P}^p(S, d)^2 : \left\lvert \int f \,d\mu - \int f \,d\nu \right\rvert+\lvert m_p(\mu)-m_p(\nu) \rvert<\varepsilon \,\right\}
\]
indexed by $f \in C_{\rm b}(S)$ and $\varepsilon>0$.
Let $\mathbf{T}^p$ denote the topology of the uniform space $(\mathbf{P}^p(S, d), \mathbf{U}^p)$.
Then $\lim_\alpha \mu_\alpha=\mu$ in $\mathbf{T}^p$ if and only if the pair of the conditions $\mu_\alpha \stackrel{\rm w}{\to} \mu$ and $m_p(\mu_\alpha) \to m_p(\mu)$ holds.

The topology $\mathbf{T}^p$ can be given in different ways.
Let $\rho$ be a distance which induces the weak topology on $\mathbf{P}(S)$.
Then the distance
\[
\rho_p(\mu, \nu):=\rho(\mu, \nu)+\lvert m_p(\mu)-m_p(\nu) \rvert,
\qquad \mu, \nu \in \mathbf{P}^p(S, d),
\]
induces the topology $\mathbf{T}^p$.
It is known as well (see Villani~\cite[Chapter~6]{Vi09} for example) that the $p$-Wasserstein distance $W_p$ on $\mathbf{P}^p(S, d)$ induces $\mathbf{T}^p$.

By Proposition~\ref{prop:cpt-open}, the topology of uniform convergence on compacta on the space $C(T; \mathbf{P}^p(S, d))$ of continuous mappings from $T$ to $\mathbf{P}^p(S, d)$ is the same if we take any one of the above constructions to give a uniformity on $\mathbf{P}^p(S, d)$.
Moreover, this space is metrizable if $T$ is locally compact $\sigma$-compact as in Corollary~\ref{cor:metrizability}~\eqref{i:luwc_metrizability}.

The sequential counterpart is given as follows:

\begin{theorem} \label{th:Wasserstein}
Let $(S, d)$ be a complete separable metric space, $T$ be a first-countable Hausdorff space, and $\rho$ be a distance inducing the weak topology on $\mathbf{P}(S)$.
Suppose that $\mu^{(n)}_t \in \mathbf{P}^p(S, d)$ for each $t \in T$ and for each $n \in \N$ and $(\mu_t)_{t \in T} \in C(T; \mathbf{P}^p(S, d))$.
Then the following are mutually equivalent:
\begin{enumerate}
\item
\eqref{eq:seq-luwc} and $\lim_{n \to \infty}\sup_{t \in K}\lvert m_p(\mu^{(n)}_t)-m_p(\mu_t) \rvert=0$ hold for every compact $K \subset T$ and for every $f \in C_{\rm b}(S)$;
\item
$\lim_{n \to \infty}\sup_{t \in K}\rho(\mu^{(n)}_t, \mu_t)=0$ and $\lim_{n \to \infty}\sup_{t \in K}\lvert m_p(\mu^{(n)}_t)-m_p(\mu_t) \rvert=0$ hold for every compact $K \subset T$;
\item
$\lim_{n \to \infty}\sup_{t \in K}W_p(\mu^{(n)}_t, \mu_t)=0$ holds for every compact $K \subset T$.
\end{enumerate}
\end{theorem}

\section{Application I: convergence of additive processes}
\label{sec:FLT}

We now observe the relation of our results to a functional limit theorem in probability theory.
Let $X=(X_t)_{t \ge 0}$ be a one-dimensional additive process, that is, an a.s.-c\'adl\'ag and stochastically continuous process with independent increments.
In the L\'evy--Khintchine representation
\begin{equation} \label{eq:L-K_1d}
\mathbb{E}[e^{i\xi X_t}]
=\exp\left[ i \gamma_t \xi + \int_{\R} \left(e^{i\xi x}-1-\frac{i\xi x}{1+x^2}\right)\frac{1+x^2}{x^2} \,\eta_t(dx) \right],
\quad \xi \in \R,
\end{equation}
let us call the family of $(\gamma_t, \eta_t) \in \R \times \mathbf{M}(\R)$, $t \ge 0$, the \emph{generating family} of $X$.
The function $t \mapsto \gamma_t$ is continuous non-decreasing with $\gamma_0=0$, and $t \mapsto \eta_t$ is setwise continuous non-decreasing, which means that $t \mapsto \eta_t(B)$ is continuous non-decreasing for every Borel set $B \in \mathscr{B}(\R)$, with $\eta_0=0$ (cf.\ \cite{It04,JS10,Sa13}).

\begin{remark}
Compared to the multidimensional formula \eqref{eq:L-K_md} below, the Gaussian coefficient in \eqref{eq:L-K_1d} is given by $\eta_t(\{0\})$; here, the value of the integrand at $x=0$ is assigned as the limit
\[
\lim_{x\to 0}\left(e^{i\xi x}-1-\frac{i\xi x}{1+x^2}\right)\frac{1+x^2}{x^2}=-\frac{\xi^2}{2}.
\]
\end{remark}

The next proposition is taken from Jacod and Shiryaev~\cite[Theorem~VII.3.4 and Corollary~VII.4.43]{JS10}, adapted to the one-dimensional case.

\begin{proposition} \label{prop:JS10}
Let $X^{(n)}=(X^{(n)}_t)_{t \ge 0}$, $n \in \N$, and $X=(X_t)_{t \ge 0}$ be one-dimensional additive processes and $(\gamma^{(n)}_t, \eta^{(n)}_t)_{t\ge 0}$ and $(\gamma_t, \eta_t)_{t\ge 0}$ be their generating families, respectively.
Then the following are mutually equivalent:
\begin{enumerate}
\item \label{JS10:J_1}
$X^{(n)} \xrightarrow{\rm d} X$ in Skorokhod's $J_1$-topology as $n \to \infty$;
\item \label{JS10:pair_pw}
The generating families satisfy
\begin{gather}
\lim_{n\to \infty}\sup_{t\in [0,T]}\lvert \gamma^{(n)}_t-\gamma_t\rvert=0
\qquad \text{for each}\ T>0\ \text{and}
\label{eq:luc_gamma} \\
\operatorname*{w-lim}_{n \to \infty}\eta^{(n)}_t=\eta_t 
\qquad \text{for each}\ t\ge 0; \label{eq:pwc_eta}
\end{gather}
\item \label{JS10:pair_lu}
The generating families satisfy \eqref{eq:luc_gamma} and
\begin{equation} \label{eq:luc_eta}
\begin{gathered}
\lim_{n \to \infty}\sup_{t \in [0,T]}\left\lvert \int_{\R} f(x) \, \eta^{(n)}_t(dx) - \int_{\R} f(x) \, \eta_t(dx) \right\rvert=0 \\
\text{for each}\ T>0\ \text{and for each}\ f \in C_{\rm b}(\R);
\end{gathered}
\end{equation}
\item \label{JS10:cf_lu}
The characteristic function $\mathbb{E}[\exp(i\xi X^{(n)}_t)]$ converges to $\mathbb{E}[\exp(i\xi X_t)]$ locally uniformly in $(\xi,t) \in \R \times [0,+\infty)$ as $n \to \infty$;
\item \label{JS10:cf_pw}
The characteristic function $\mathbb{E}[\exp(i\xi X^{(n)}_t)]$ converges to $\mathbb{E}[\exp(i\xi X_t)]$ locally uniformly in $t \in [0,+\infty)$ as $n \to \infty$ for each $\xi \in \R$.
\end{enumerate}
\end{proposition}

In Proposition~\ref{prop:JS10}, let $\mu^{(n)}_t$ and $\mu_t$ be the distributions of $X^{(n)}_t$ and $X_t$, respectively.
Then Condition \eqref{JS10:cf_lu} or \eqref{JS10:cf_pw} is equivalent to $(\mu^{(n)}_t)_{t \ge 0} \to (\mu_t)_{t \ge 0}$ in $C([0,+\infty); \mathbf{P}(\R))$ by Theorem~\ref{th:Levy_continuity}.
Moreover, \eqref{eq:luc_eta} means that $(\eta^{(n)}_t)_{t \ge 0} \to (\eta_t)_{t \ge 0}$ in $C([0,+\infty); \mathbf{M}(\R))$.
Hence, if a sequence of additive processes converges to an additive process, then the finite measures in the L\'evy--Khintchine representation as well as the marginal distributions of the processes converge locally uniformly in time parameter.

\begin{remark} \label{rem:Polya}
In Proposition~\ref{prop:JS10}, the conditions \eqref{eq:pwc_eta} and \eqref{eq:luc_eta} are equivalent because $\eta^{(n)}_t$ and $\eta_t$ are setwise continuous and non-decreasing in $t$.
In fact, a theorem of P\'olya~\cite{Po20} asserts that, if a sequence of right-continuous non-decreasing functions $F_n$, $n \in \N$, on $[0,+\infty)$ converges to a continuous function $F$ pointwise with $F_n(0)=F(0)=0$, then $(F_n)_{n \in \N}$ converges to $F$ locally uniformly on $[0,+\infty)$.
From this theorem we easily see that \eqref{eq:pwc_eta} implies \eqref{eq:luc_eta}, decomposing a test function $f \in C_{\rm b}(\R)$ into the positive and negative parts.
Since the marginal distributions $\mu^{(n)}_t$ of the processes $X^{(n)}_t$ have no monotonicity in $t$, we cannot drop $t$-uniformity in the convergence of $((\mu^{(n)}_t)_{t \ge 0})_{n \in \N}$.
\end{remark}

In relation to Remark~\ref{rem:Polya}, there is still another condition equivalent to \eqref{eq:pwc_eta} or \eqref{eq:luc_eta}.
To state this, we recall that $t \mapsto \eta_t(B)$ is continuous, non-decreasing for every $B \in \mathscr{B}(\R)$.
Then there exists a unique $\sigma$-finite Borel measure $\mathcal{H}$ on $\R \times [0,+\infty)$ that enjoys
\begin{align}
&\mathcal{H}(B \times (s,t]) = \eta_t(B) - \eta_s(B)
&&\text{for}\ B \in \mathscr{B}(\R) \ \text{and}\ 0 \le s \le t, \label{eq:two-dim_DF} \\
&\mathcal{H}(\R \times [0,t])<+\infty \ \text{and}\ \mathcal{H}(\R \times \{t\})=0
&&\text{for}\ t \ge 0. \label{eq:t-atomless}
\end{align}
Indeed, $(x, t) \mapsto \eta_t((-\infty, x]$ ($=\mathcal{H}((-\infty, x] \times [0,t])$) is a two-dimensional distribution function (cf.\ Billingsley~\cite[Theorem~12.5]{Bi12}).
Using this $\mathcal{H}$, we show the following:

\begin{theorem}
Let $(\eta_t)_{t \ge 0}$ and $(\eta^{(n)}_t)_{t \ge 0}$, $n \in \N$, be mappings from $[0,+\infty)$ to $\mathbf{M}(\R)$ with $\eta_0=\eta^{(n)}_0=0$ which are nondecreasing and continuous in $t$ setwise.
Let $\mathcal{H}$ be the $\sigma$-finite measure defined by \eqref{eq:two-dim_DF} and \eqref{eq:t-atomless}.
For each $n$, let $\mathcal{H}^{(n)}$ be the measure defined by \eqref{eq:two-dim_DF} and \eqref{eq:t-atomless} with the superscript ${}^{(n)}$ put on $\eta$ and $\mathcal{H}$.
Then the following are equivalent:
\begin{enumerate}
\item \label{i:eta}
\eqref{eq:pwc_eta} or \eqref{eq:luc_eta} holds;
\item \label{i:curl_H}
$\mathcal{H}^{(n)}\rvert_{\R \times [0,T]} \xrightarrow{\rm w} \mathcal{H}\rvert_{\R \times [0,T]}$ holds for each $T>0$.
\end{enumerate}
\end{theorem}

\begin{proof}
Since
\[
\int_{\R \times [0,T]} f(x) \,\mathcal{H}^{(n)}(dx\,dt) = \int_{\R} f(x) \,\eta^{(n)}_T(dx),
\quad f\in C_{\rm b}(\R),
\]
Condition \eqref{i:curl_H} implies \eqref{eq:pwc_eta}.
Conversely, suppose \eqref{eq:pwc_eta}.
Then repeating the proof of $\rm (G2)\Rightarrow (G3)$ in Hasebe and Hotta~\cite[Proposition~2.6]{HH22}, we see that $\mathcal{H}^{(n)}\rvert_{\R \times [0,T]}$ converges vaguely to $\mathcal{H}\rvert_{\R \times [0,T]}$ as $n \to \infty$.
Moreover, the total mass $\mathcal{H}^{(n)}(\R \times [0,T])=\eta^{(n)}_T(\R)$ converges to $\mathcal{H}(\R \times [0,T])=\eta_T(\R)$ by \eqref{eq:pwc_eta}.
Hence \eqref{i:curl_H} follows.
\end{proof}

For multidimensional additive processes, the same property as Proposition~\ref{prop:JS10} holds true, except that the pair $(\gamma_t, \eta_t)$ in \eqref{eq:L-K_1d} should be replaced by the characteristic triplet.
Namely, the following expression is adopted:
\begin{equation} \label{eq:L-K_md}
\mathbb{E}[e^{i\xi \cdot X_t}]
=\exp\left[ i \gamma_t \cdot \xi - \frac{1}{2}\xi \cdot C_t\xi + \int_{\R} \left(e^{i\xi x}-1-i\xi \cdot h(x) \right) \,\nu_t(dx) \right],
\quad \xi \in \R^d.
\end{equation}
Here, $C_t$ is the Gaussian covariance matrix, $h \colon \R^d \to \R^d$ is a truncation function~\cite[p.75]{JS10}, and $\nu_t$ is a L\'evy measure.
Then
\[
\tilde{C}_t := C_t - \int_{\R^d} h(x) \otimes h(x) \,\nu_t(dx)
\]
is called the modified second characteristic~\cite[p.115]{JS10}.
The multidimensional version of Proposition~\ref{prop:JS10} states that the locally uniform convergences of additive processes, of modified triplets, and of marginal distributions (of the processes) are all equivalent; for the precise statement, see Jacod and Shiryaev~\cite[Theorem~VII.3.4 and Corollary~VII.4.43]{JS10}.

\begin{remark}
An analogue to Proposition~\ref{prop:JS10} is formulated for a sequence of random walks constructed by an infinitesimal triangular array in Hata's thesis~\cite{Ha23}.
Similar results were proved even on a Lie group; see Feinsilver~\cite{Fe78}.
\end{remark}

\section{Application II: Bercovici--Pata bijection for non-commutative additive processes}
\label{sec:BP_bij}

In this section, we describe the role of our results in non-commutative probability theory, which has motivated us originally.

In non-commutative probability, we have five standard ``independences'': classical, free, boolean, monotone, and anti-monotone ones.
The symbols $\ast$, $\boxplus$, $\uplus$, $\rhd$, and $\lhd$ denote convolution of Borel probability measures on $\R$ with respect to these independences, respectively.
Bercovici and Pata~\cite{BP99} gave a canonical bijection between the $\ast$-, $\boxplus$-, and $\uplus$-infinitely divisible distributions, which was later extended to the $\rhd$-infinitely divisible distributions \cite{AW14}.

Hasebe and Hotta~\cite{HH22} gave a dynamical version of the Bercovici--Pata bijection, which is a correspondence between the marginal distributions of classical, free, boolean, and monotone additive%
\footnote{As the group $\mathbb{T}=\{\, e^{i\theta} : \theta \in \R \,\}$ is concerned, they are actually ``multiplicative''.}
processes \emph{on the unit circle $\mathbb{T}$}.
Moreover, they proved that the bijection is continuous with respect to the convergence \eqref{eq:seq-luwc} with $S=\mathbb{T}$ and $T=[0,+\infty)^2_{\le}:=\{\, (s,t) : 0 \le s \le t \,\}$.
In our companion paper~\cite{HHM23}, we shall establish similar results for processes on $\R$.
Below we announce part of the results related to the present article.

Let $\star$ denote one of the convolutions $\ast$, $\boxplus$, $\uplus$, and $\rhd$.
A family $(\mu_{s,t})_{0 \le s \le t} \in C([0,+\infty)^2_{\le}; \mathbf{P}(\R))$ is called a \emph{$\star$-convolution hemigroup} if $\mu_{s,s}=\delta_0$, $s \ge 0$, and if $\mu_{s,u}=\mu_{s,t} \star \mu_{t,u}$, $0 \le s \le t \le u$.
If $(X_t)_{t \ge 0}$ is a $\star$-additive process, then the distributions $\mu_{s,t}$ of the increments $X_t-X_s$, $0 \le s \le t$, form a $\star$-convolution hemigroup.
We define $\mathrm{CH}(\star)$ as the set of $\star$-convolution hemigroups and set $\mathrm{CH}_2(\star):=\{\, (\mu_{s,t})_{0 \le s \le t} \in \mathrm{CH}(\star) : \mu_{s,t} \in \mathbf{P}^2(\R) \,\}$.
In the next proposition, $\mathbf{P}^2(\R)=\mathbf{P}^2(\R, d)$ is equipped with the topology $\mathbf{T}^2$ in Section~\ref{sec:vague_moment} with $d$ being the Euclidean distance.

\begin{proposition}[{\cite{HHM23}}]
\label{prop:forthcoming}
There exists a ``canonical'' bijection between $\mathrm{CH}_2(\ast)$ and $\mathrm{CH}_2(\rhd)$, which is homeomorphic with respect to the relative topologies in $C([0,+\infty)^2_{\le}; \mathbf{P}^2(\R))$.
Similarly, there exists a bijection between $\mathrm{CH}(\ast)$, $\mathrm{CH}(\boxplus)$, and $\mathrm{CH}(\uplus)$, homeomorphic with respect to the relative topologies in $C([0,+\infty)^2_{\le}; \mathbf{P}(\R))$.
\end{proposition}

Basically, the bijection $\mathrm{CH}_2(\ast) \ni (\mathfrak{c}_{s,t})_{s,t} \mapsto (\mathfrak{m}_{s,t})_{s,t} \in \mathrm{CH}_2(\rhd)$ in Proposition~\ref{prop:forthcoming} is obtained by identifying the generating family of the classical additive process whose marginal distributions are $\mathfrak{c}_{0,t}$, $t \ge 0$, with the ``generating family'' of $(\mathfrak{m}_{s,t})_{s,t}$, which will be newly defined in the paper~\cite{HHM23}.
Then the choice of the topology in Proposition~\ref{prop:forthcoming} should be natural in view of Proposition~\ref{prop:JS10}.

\begin{remark}
Hasebe and Hotta~\cite{HH22} showed a version of Proposition~\ref{prop:forthcoming} on $\mathbb{T}$ in the following form:
there exists a certain bijection $\beta$, and it is \emph{sequentially bi-continuous}, i.e.,
\[
(\mu^{(n)}_{s,t})_{s,t} \to (\mu_{s,t})_{s,t} \ (n \to \infty)
\iff
\beta\bigl((\mu^{(n)}_{s,t})_{s,t}\bigr) \to \beta\bigl((\mu_{s,t})_{s,t}\bigr) \ (n \to \infty).
\]
By virtue of Corollary~\ref{cor:metrizability}~\eqref{i:luwc_metrizability}, this now implies genuine continuity.
\end{remark}

\section{Other related studies}

We have regarded the domain of definition $T$ of measure-valued functions as the set of time parameters, but another interpretation is possible.
For example, Kawabe~\cite{Ka99} interpreted a continuous function $\lambda \colon T\to \mathbf{M}(S)$, $t \mapsto \lambda(t, ds)$ as a transition probability kernel from $T$ to $S$.
Let $(\mu^\alpha)_{\alpha \in A}$ be a net of initial distributions on $T$ and $(\lambda^\alpha)_{\alpha \in A}$ be a net of transition probability kernels.
He studied, using uniform structures on $\mathbf{M}(S)$ and their relation to the topology of $C(T; \mathbf{M}(S))$, in what case the net of the compound probability measures $\int_T \lambda^\alpha(t, \cdot) \,\mu^\alpha(dt)$ converge.

Bengs and Holzmann~\cite{BH19a,BH19b} studied the uniform convergence of a sequence of families of distributions $(\, \mu^\vartheta_n : \vartheta\in\varTheta \,)$, $n\in \mathbb{N}$, on $\mathbb{R}^d$ in the context of statistics.
They assumed that the parameter set $\varTheta$ has no topological structure and the limit family $(\, \mu^\vartheta : \vartheta\in\varTheta \,)$ is ``uniformly absolutely continuous'' with respect to a fixed atomless distribution.
Their assumption is thus completely different from ours; nevertheless, they proved the corresponding version of L\'evy's continuity theorem as well as the portmanteau theorem and central limit theorem.

\section*{Acknowledgments}

The authors would like to thank the anonymous referees for their valuable comments on old versions of this paper.
This work is supported by JSPS KAKENHI Grant Numbers  18H01115, 19K14546, 20H01807, 20K03632, 22K20341, 23K25775, 24K16935, and JSPS Open Partnership Joint Research Projects grant no.~JPJSBP120209921.


\end{document}